%This is an ams-tex paper using style files for the Proceedings of the AMS  
 
\input amstex
\loadmsbm
\loadeufm
\documentstyle{amsppt}

\def\R{{\Bbb R}}   % real numbers
\def\Z{{\Bbb Z}}   % integers
   % rational numbers
\def\C{{\Bbb C}}   % complex numbers
\def\H{{\Bbb H}}   % quaternionic numbers 

\def\Map{{\roman{Map}}}

\def \gr#1{ G_{#1}(V^{\infty})}
\def \g#1{G_{#1} (\H ^{\infty})}

\def\A{{\Cal A}} 
\def\AH{{\Cal A}_{U(1)}}
\def\B{{\Cal B}} 
 
\def\BH{{\Cal B}_{U(1)}}
\def\BHO{{{\Cal B}_{U(1)}^0}} 
\def\BO{{\Cal B}^0}
\def\G{{\Cal G}}
\def\GH{{{\Cal G}_{U(1)}}} 
\def\GHO{{{\Cal G}_{U(1)}^0}} 
\def\GO{{{\Cal G}^0}}

\topmatter
\title The Homotopy Type of hyperbolic Monopole Orbit Spaces\endtitle
\author Ursula Gritsch\endauthor
\leftheadtext{}%
\rightheadtext{}%

\address DPMMS, University of Cambridge, 16 Mill Lane, Cambridge, CB2 1SB, 
U.K.\endaddress
\email ursula\@dpmms.cam.ac.uk\endemail


\abstract We prove that the space $\BHO$ of equivalence classes of
$U(1)$-invariant connections on some $SU(2)$-principle bundles
over $S^4$ is weakly homotopy equivalent to a component of the 
second loop space $\Omega^2 (S^2)$.\endabstract

\thanks 
This note is part of my PhD-thesis written at Stanford University,
1997.  I thank my advisor Ralph Cohen for constant support and
encouragement and the Studienstifung des deutschen Volkes for a
dissertation fellowship.\endthanks

\date \enddate


\endtopmatter

\document

\head 1. Introduction
\endhead 

In this note we prove that the space $\BHO$ of gauge equivalence
classes of $U(1)$-invariant connections on an $SU(2)$-principle bundle
over the 4-sphere $S^4$ of even Chern class is weakly homotopy
equivalent to one component of the second loop space of the sphere
$S^2$. In other words we prove:
\proclaim {Theorem 3.2.1}
$$
\BO_{U(1)} \simeq_{\roman w} \Omega^2_k S^2 \; . 
$$
\endproclaim

This $U(1)$-action on $S^4 \subset \R^5$ is just the rotation in the 
$\langle x_1 , x_2 \rangle$-plane leaving the $x_3$, $x_4$ and $x_5$-coordinates fixed. 

We use an equivariant version (see [G, theorem 2.3.1]) of a theorem in
[AB, prop. 240, p.540 ] or [DK, prop. 5.1.4, p.174], which identifies
the invariant orbit space $\BHO$ as an equivariant mapping space. We
then analyze this mapping space using methods from equivariant
homotopy theory. Usually theorems of this kind are proved using tools
from analysis and gauge theory.

Anti self dual connections invariant under this $U(1)$-action on the
4-sphere $S^4$ are called hyperbolic monopoles and were first
studied by Atiyah in [A].  For this reason we call the space $\BHO$
the hyperbolic monopole orbit space.

In [A, section 2, p.3] Atiyah shows that anti self dual
connections on the 4-sphere $S^4$ invariant under this particular
$U(1)$-action are in one-to-one correspondence with pairs $(A, \Phi )$, 
where $A$ is a connection on the trivial $SU(2)$-bundle over the
hyperbolic 3-space ${\Cal H}^3$, the Higgs field $\Phi$ is a section
of the adjoint bundle ${\Cal H}^3 \times \goth{su}(2)
\rightarrow {\Cal H}^3$, and the pair $(A, \Phi )$ satisfies the 
Bogomolny equation  
$D\Phi = \ast F_A$ on the space ${\Cal H}^3$ and has finite energy. 

In [T1] Taubes investigates the Euclidean analogue. He studies
pairs $(A, \Phi )$ where $A$ is a connection on the trivial
$SU(2)$-bundle over $\R^3$, the Higgs field $\Phi$ is a section of the
adjoint bundle $\R^3 \times \goth{su}(2) \rightarrow \R^3$ satisfying a
regularity condition and the pair $(A,\Phi )$ has finite energy with
respect to the Yang-Mills-Higgs functional.  In [T1, Theorem A1.1,
p.347] he shows that the configuration space $\B$ (that is the space
of smooth pairs $(A, \Phi )$, where $A$ and $\Phi$ are as above)
modulo the action of the based gauge group of the trivial
$SU(2)$-bundle over $\R^3$ is homotopy equivalent to the space of
smooth maps $f: S^2 \rightarrow S^2$.  The space $\B$ fibers over the
2-sphere $S^2$ and C.H. Taubes shows that the fiber $\hat{\B}$ over
the north pole is homotopy equivalent to the space of based smooth
maps $f: S^2 \rightarrow S^2$.

\head 2. The hyperbolic monopole orbit space $\BHO$
\endhead

\subhead {\rm 2.1} Notation and terminology
\endsubhead 

Let $S^4$ denote the 4-sphere in $\R^5$ and let $ \eta = (P
\rightarrow S^4)$ be an $SU(2)$-principal bundle over the manifold
$S^4$. Denote by $\A$ the space of smooth connections on the bundle
$\eta$. Let $\G$ be the group of smooth gauge transformations of this
bundle. This means that elements $g \in \G$ are smooth automorphisms
of $\eta$. Let $\GO$ be the subgroup of the group $\G$ of
gauge transformations whose elements are the identity over a given
base point $m \in S^4$.

The circle $U(1)$ acts on the 4-sphere 
$S^4 \subset \R^5$ by the
rotation in the $\langle x_1 ,x_2 \rangle$-plane leaving the $x_3$, $x_4$ and
$x_5$-coordinates fixed.
We assume that the $U(1)$-action on the base $S^4$
lifts to the total space $P$ of the bundle $\eta$.  It follows from
[Ba, proposition 2.1, p.434] and [A, 2.6, p.6] that this is the case
if and only if 
the second Chern class
$c_2 (\eta ) = n$ is even. In [Ba, p.430] P. Braam shows that the
induced action of the non trivial double cover ${\tilde U(1)}$ of
$U(1)$ on the 4-sphere $S^4$ always lifts to the bundle $\eta$
regardless of the parity of $n$.

By abuse of notation we denote by $ \eta = ( E= P \times_{SU(2)} \C^2
\rightarrow S^4 )$ the two dimensional complex vector bundle associated 
to the principle bundle via the standard representation of $SU(2)$ 
on the complex vector space $\C^2$. The complex vector 
bundle $\eta$ can be given the structure of a quaternionic line 
bundle using the isomorphisms $SU(2) \cong Sp(1) $ and 
$ \C^2 \cong \H$ as $SU(2) \cong Sp(1)$ left modules. 
(We consider the quaternions $\H$ as a right quaternionic vector space).

We define the invariant gauge group $\GH$  
to be the subgroup of the gauge group $\G$ such that each gauge 
transformation commutes with the action of the Lie group $U(1)$ 
on the total space $P$ of the bundle $\eta$. 
Analogously we define the based invariant gauge group 
$\GHO$ to be the equivariant gauge transformations which are the 
identity over a given based orbit $ {\Cal O}_m \subset S^4$. 

The groups $\G$ and $\GO$ act naturally on the space of connections 
$\A$ from the right by the pull back of connections. 
We define the orbit space 
$ \B $ to be the space $ \A / \G $ of connections modulo the action of the 
gauge group and the based orbit space $\BO $ to be the space $ \A / \GO $  
of connections modulo the action of the based gauge group. 

The left action of the group $U(1)$ on the bundle $\eta$ induces 
a right action on the space of connections $\A$ again by the pull back 
of connections. The fixed points of this action are called 
$U(1)$-invariant connections. We denote the space of fixed points 
by $\AH$. The invariant groups $\GH$ and $\GHO$ act on the space 
of invariant connections $\AH$ as above. As in the non-equivariant 
setting we define the invariant orbit space 
$ \BH $ to be the space $ \AH / \GH $ of invariant connections modulo the 
invariant gauge group and the based invariant orbit space 
(the hyperbolic monopole orbit space)
$ \BHO $ to be the space $ \AH / \GHO $ of invariant connections 
modulo the based invariant gauge group.

It is customary to complete spaces of connections in the 
$L^{2,2}$-Sobolev norm
and elements of the gauge group in the $L^{3,2}$-Sobolev norm 
using a fixed connection on the bundle $\eta$.  
However for the purpose of this paper this is not necessary. So the 
reader can either think of smooth connections and smooth gauge 
transformations or of $L^{2,2}$-connections and $L^{3,2}$-gauge 
transformations. 

\subhead {\rm 2.2} The mapping space model of the hyperbolic monopole orbit space $\BHO$\endsubhead  

In this section we briefly describe the model for the hyperbolic monopole 
orbit space $\BHO$ in terms of an equivariant mapping space given in 
[G, theorem 2.3.1]. 
It is a direct generalization of the model for the 
orbit space $\BO$ given in [DK] and [AB] to the invariant situation. 
The results in this section hold for any action of a compact Lie group $H$ 
on a principle bundle $\eta =(P \rightarrow M)$ over any 
compact $H$-manifold $M$. However we state them only for $SU(2) \cong Sp(1)$-principle 
bundles over $S^4$ and $H=U(1)$. 

We first introduce a few notions from equivariant bundle theory.  Let
$V$ be a fixed quaternionic $U(1)$-representation.  That means $V$ is
a quaternionic right vector space with $U(1)$ left action.  Let $
G_1(V^k)$ be the Grassmannian of quaternionic lines in the vector
space $V^k = \underbrace{V \oplus \dots \oplus V}_{k \ \roman times} $
for $k >1$. The group $U(1)$ acts naturally on the space $ G_1(V^k)$
because it acts on the space $V^k$.  ($ G_1(V^k)$ is just a
quaternionic projective space with a $U(1)$ left action).
  
We denote by 
$G_1(V^{\infty}) = \lim_{k \rightarrow \infty } G_1(V^k)$ the direct limit as 
an $U(1)$-equivariant space.  
Let $ \gamma_1 (V^k ) = ( E_1 (V^k) \rightarrow G_1 (V^k ) )$ 
denote the canonical unitary bundle over the 
Grassmannian. The fiber over a line $P \in G_1 (V^k )$ 
are the points $ p \in P$. 
The $U(1)$-action on the Grassmannian lifts naturally to an $U(1)$-action on the 
total space of the canonical bundle giving this bundle the structure of an 
$U(1)$-equivariant bundle. We take the limit 
$ E_1(V^{\infty}) = \lim_{k \rightarrow \infty } E_1(V^k)$
and obtain the $U(1)$-equivariant bundle 
$ \gamma_1 (V^{\infty}) = ( \pi :E_1(V^{\infty})  \rightarrow 
G_1(V^{\infty}) )$. 

\definition{Definition 2.2.1} 
We call an $U(1)$-equivariant quaternionic vector bundle 
$\eta = (E \rightarrow S^4 )$ subordinate to the 
representation $V$ if 
for every $m \in S^4 $ the isotropy representation 
of the isotropy group $U(1)_m$ on the fiber over $m$ is equivalent to 
a sub-module of the  
$U(1)_m$-module $V$ induced by the given 
$U(1)$-module $V$.\enddefinition 
The bundle 
$ \gamma_1 (V^{\infty}) = ( \pi :E_1(V^{\infty}  \rightarrow 
G_1(V^{\infty}) )$ classifies 
vector bundles subordinate to the representation $V$. 
This is proved by Wasserman in [Wa, section 2, p.132] for the case of 
real vector bundles and by Segal in [Se, section 1, p.131] for complex 
vector bundles. These proofs also hold in the quaternionic category. 

Similarly if $\eta = (P \rightarrow M)$ is an $U(1)$-equivariant 
$SU(2) \cong Sp(1) $-principle bundle such that the associated vector bundle 
$\eta = (E = P \times_{Sp(1)}\H \rightarrow S^4)$ is subordinate to the 
representation $V$ then the frame bundle associated to the bundle 
$ \gamma_1 (V^{\infty}) = ( \pi : E_1(V^{\infty})  
\rightarrow G_1(V^{\infty}) )$ 
classifies the principle bundle $\eta$. 
\definition{Definition 2.2.2}
We call 
$\gr1 $ a classifying space for the bundle $\eta$ provided it 
is subordinate to the $U(1)$-representation $V$.\enddefinition

\remark{Remark \rom{2.2.3}} Fix a point $e \in S^4$ which is fixed under the $U(1)$-action.
Also we fix a trivialization of the bundle $\eta$ over the point $e$
once and for all. Then there is a point $\ast \in G_1 (V^{\infty})$
such that the bundle $(\eta , \eta / e )$ is classified by the pointed
space $ (G_1 (V^{\infty}) , \ast)$.  Of course the base point $ \ast
\in G_1 (V^{\infty})$ also has to be a fixed point of the circle
action on $G_1 (V^{\infty})$.\endremark

\proclaim {Proposition 2.2.4} Assume that the $U(1)$-equivariant 
$SU(2)$-principle bundle $\eta = (P \rightarrow S^4)$ 
is subordinate to the $U(1)$-module $V$. Then 
there is a weak homotopy equivalence 
$$  
\BHO \simeq_{\roman w} 
\Map _{U(1)}^0 (S^4 , \; \gr1))^{\eta} 
$$ 
where the right hand side denotes the component of the mapping space of maps 
which classify the bundle $(\eta , \eta / e )$. 
(The mapping space is given the compact-open topology). 
\endproclaim

The proof is given in [G, theorem 2.3.1]. It is a direct
generalization to the invariant situation of the corresponding theorem
in [DK, prop. 514, p.174].

\head 3. The weak homotopy type of the hyperbolic monopole orbit 
space $\BHO$\endhead 

\subhead {\rm 3.1} Background and definitions\endsubhead 

In this section we apply proposition 2.2.4 to compute the weak
homotopy type of the hyperbolic monopole orbit space $\BHO$.  Recall
that the $U(1)$-action on $S^4 \subset \R^5$ is just the rotation in
the $\langle x_1 ,x_2 \rangle$-plane, leaving the $x_3 , x_4 , $ and $x_5$-coordinates 
fixed.  The fixed point set is the two-sphere consisting of the points
with $x_1 =x_2 =0$ and hence $x_3^2 + x_4^2 + x_5^2 = 1$.

However, to use methods from equivariant homotopy theory it is better to 
think of the sphere $S^4$ as the second suspension of a certain 
representation sphere. 
For an orthogonal representation $V$ of $U(1)$ 
let $S(V)$ denote the unit sphere in the representation $V$ 
and $S^V$ the 1-point compactification. Both 
spaces carry an induced $U(1)$-action. The 2-sphere with the standard 
circle action given by rotation around the axis through the north and south 
pole is the unit sphere $S^c = S(1+c)$. Here $c$ denotes the standard 
representation of $U(1)$ on $\C$ and $1$ denotes the trivial 1-dimensional  
representation. The 2-sphere with the trivial circle action is the unit sphere 
$S^2 = S(1+1+1)$. 

For any space 
$X$ let $ S^0 \ast X$ denote the (unreduced) suspension of $X$. Hence 
$S^0 \ast X = {X \times {\roman I} \over 
X \times \{ 0 \} , X \times \{ 1 \} }$  
where the space $I$ denotes the unit interval $[0,1]$ and we collapse the 
subspaces $ X \times \{ 0 \} $ and $  X \times \{ 1 \}$ 
to two different points. 
The space $S^1 \ast X$ denotes the two fold suspension of $X$, i.e. 
we suspend the space $X$ twice. 

If the space $X$ is acted on by the group $U(1)$ then the two-fold 
suspension $S^1 \ast X$ is acted on naturally by $U(1)$ as well, namely: 
$ \lambda \cdot [s,t,v] = [s,t, \lambda \cdot v]$ for elements 
$s$, $t \in I$,  
$v \in S^2$ and $\lambda \in U(1)$.    
There is a well-known homeomorphism between the 4-sphere $S^4$ and the second 
unreduced suspension of the 2-sphere. Choosing this 2-sphere to be $S^c$ 
induces a $U(1)$-action on $S^4$, i.e. $S^4 = S^1 \ast S^c$ as a 
$U(1)$-space.  

\remark{Remark \rom{3.1.1}}
We could also regard $S^4$ as the third unreduced suspension of the
circle $U(1)$ where $U(1)$ acts on itself by multiplication. However
this is not very useful when using proposition 2.2.4: The action of
$U(1)$ on itself has no fixed points but the mapping space in
proposition 2.2.4 is a based space.  Since methods from equivariant
homotopy theory when applied to based mapping spaces work well when
one has fixed points we prefer to work with the suspension of a
$U(1)$-space with fixed points; for example the second suspension of
the 2-sphere $S^c$.  (This is used in the paragraph following
3.2.3).\endremark

The $U(1)$-action on $S^c$ fixes the north and south poles $N$ and
$S$, and hence the circle action on $S^4$ fixes the space $ S^1 \ast (
N \cup S) \cong S^2$.  Therefore the circle action on $S^4$ has a 2-sphere
inside the 4-sphere as the fixed point set.
 
In order to apply proposition 2.2.4 we need to compute the family of
isotropy representations of the $U(1)$-action on the bundle $\eta $
and to determine a classifying space for the bundle $\eta$.  The
following analysis can be found in [A, p.4 and 5].

There are two types of isotropy groups. If a point $ m \in S^4$ lies
in the complement of the fixed point set then the isotropy group
$U(1)_m$ is trivial and hence the corresponding isotropy
representation is trivial. If a point $m \in S^4$ lies in the fixed
point set $S^2 \subset S^4$ the isotropy group is the whole circle
$U(1)$. 

Let $E_m$ denote the fiber of the two dimensional complex
vector bundle $\eta$ over the fixed point $m$. Since the structure
group of $\eta$ is $SU(2)$ the $U(1)$-action on the vector space $E_m$
is conjugate to the standard two dimensional complex representation of
$U(1)$ given by the matrix
$$ 
\pmatrix \lambda^p & 0 \\ 
              0 & \lambda^{-p} \\
\endpmatrix
$$ 
for some non negative integer $p$ called the mass. 

Let $A$ be any $U(1)$-invariant connection on the bundle $\eta$. 
For two different points $m_1$ and $m_2$ in the fixed point set $S^2$ 
the parallel transport induced by the connection $A$ defines an 
isomorphism as $U(1)$-modules of the fibers $E_{m_1}$ and $E_{m_2}$. 
Hence 
the mass $p$ is independent of the choice of the point $m$ in 
the fixed point set $S^2 \subset S^4$. 

The bundle $\eta $ restricted to the fixed point set $S^2$ splits 
into the direct sum of two complex line bundles $l_p$ and $l_{-p}$.  
The fiber of the bundle $l_p$ at a point $m \in S^2$ consists of 
the eigenspace with eigenvalue $\lambda^p$ for $\lambda \in U(1)$.  
Likewise the fiber of the bundle $l_{-p}$ consists of eigenspaces 
with eigenvalues $\lambda^{-p}$ for $\lambda \in U(1)$. 
The line bundles  $l_p$ and $l_{-p}$ are complex conjugate to 
each other. 
\definition{Definition 3.1.1}
Let $k$ be the first Chern class $c_1 (l_p )$ of the bundle 
$l_p$. The integer $k$ is called the charge.\enddefinition 
It is well known [A, 2.6, p.6] that there is the equality  
$ n = 2pk$ between the mass, charge and the second Chern 
class of the bundle $\eta$. 

If we consider the bundle $ \eta$ as a quaternionic line bundle then 
over the fixed point set $S^2$ the isotropy representation of $U(1)$ 
is conjugate to the one dimensional quaternionic $U(1)$-module $\H$ where 
an element $\lambda \in U(1)$ acts by left multiplication with the 
element $\lambda^p$ on the quaternions $\H$. We say that as an $U(1)$-module 
$\H$ has weight $p$. 
Hence as a quaternionic line bundle the bundle $\eta$ is subordinate 
to the $U(1)$-module $\H$ of weight $p$ and therefore 
the equivariant Grassmannian 
$ \g {1} = \H P^{\infty}$ classifies the bundle $\eta$. 

We choose the point $e := [1,1,N] \in S^1 \ast S^c$ as the based orbit where 
$N \in S^c$ is the north pole. 
Also we fix a trivialization of the bundle $\eta$ over the point $e$
once and 
for all. Then, as mentioned in 2.2.3, 
there is a point $\ast \in \g {1}$ such 
that the bundle $(\eta , \eta / e )$ is classified by 
the pointed space $ (\g {1} , \ast)$ for some point 
$\ast \in \g {1}$. Of course the base point 
$ \ast \in \g {1}$ also has to be a fixed point of the circle action 
on $ \g {1}$.  

\subhead {\rm 3.2} Proof of the main theorem\endsubhead 

For any topological space $X$ denote by 
$
\Omega^2 X = \Map ^0 ( S^2 , X) 
$ 
the based continuous maps from the 2-sphere $S^2$ into $X$ 
together with the compact-open topology. This space has components labeled 
by the second homotopy group $\pi_2 (X)$ of $X$. For some element 
$k \in \pi_2 (X)$ let $ \Map ^0_k ( S^2 , X) $ denote the component 
labeled by $k$. 
Recall that 
the hyperbolic monopole orbit space 
$ \BO_{U(1)} = \A _{U(1)} /\G^0_{U(1)} $ 
is the space of $U(1)$-invariant connections on the bundle $\eta$ modulo 
the based equivariant gauge group $\G^0_{U(1)}$. 

\proclaim {Theorem 3.2.1}
There is a weak homotopy equivalence 
$
\BO_{U(1)} \simeq_{\roman w} \Omega^2_k S^2 ,  
$ 
where $ k$ is the charge of the bundle $\eta$ as defined in section 3.1 
and $  \Omega^2_k S^2$ denotes the component of the second loop space 
$ \Omega^2 S^2 $ labeled by the charge $k$. 
\endproclaim

\demo{Proof} Proposition 2.2.4 already gives a weak homotopy equivalence 
$$ 
\BO_{U(1)} \simeq_{\roman w} \Map ^0_{U(1)} (S^4 , \, \g {1} )^{\eta} .
\eqno{(3.2.2)} 
$$ 
Here the superscript $\eta$ in the mapping space denotes the component 
of maps which classify the bundle $\eta$. 
We will first show that there is a homotopy equivalence 
$$ 
 \Map ^0_{U(1)} (S^4 , \, \g {1}) \simeq \Omega^2 S^2 . 
\eqno{(3.2.3)}
$$ 
In section 3.3 we will identify which components correspond to each 
other under this homotopy equivalence.\enddemo 

For two pointed topological spaces $ (X,x_0 ) $ and $(Y,y_0 )$ define
the smash $ X \wedge Y$ to be the pointed space ${X \times Y \over X
\times \{ y_0 \} \cup \{ x_0 \} \cup Y }$.  If the spaces $X$ and $Y$
are acted on by the circle group $U(1)$ with fixed points $x_0$ and
$y_0$ then the smash $X \wedge Y$ is acted on in a natural fashion by
the circle $U(1)$ as well. Hence the circle group $U(1)$ acts on the
smash $S^2 \wedge S^c$ where we have chosen the north pole in both
2-spheres as the base point. Here we recall that $S^c$ denotes the
2-sphere with the standard circle action and $S^2$ the 2-sphere with
the trivial circle action.

We have identified $S^4 = S^1 \ast S^c$ and there is a well-known
$U(1)$-equivariant homotopy equivalence $S^1 \ast S^c \simeq S^2
\wedge S^c$. Hence we have to analyze the mapping space $ \Map
^0_{U(1)} (S^2 \wedge S^c , \, \g {1})$.  Taking $U(1)$-fixed points
of the equivariant homeomorphism
 $$
\split 
 \Map ^0 (S^2 \wedge S^c , \, \g {1}) 
&\cong \Map ^0 (S^2 , \, \Map^0 ( S^c , \, \g {1} )) \\ 
&\cong  \Omega ^2 \Map^0 ( S^c , \, \g {1} )\endsplit
$$ 
gives the homotopy equivalence 
$$ 
\Map ^0_{U(1)} (S^4 , \, \g {1} ) \simeq 
\Omega ^2 \Map^0_{U(1)} ( S^c , \, \g {1} ) \; . 
\eqno{(3.2.4)}
$$ 
The proof of the homotopy equivalence 3.2.3 follows now from 
lemma 3.2.5 below. 

\proclaim {Lemma 3.2.5}
There is a homotopy equivalence 
$
 \Map^0_{U(1)} (S^c , \, \g {1} ) \simeq S^2 \; . 
$  
\endproclaim

\demo{Proof}
The sphere $S^c$ is $U(1)$-equivariantly obtained by gluing in 
along the north-and 
south pole a $U(1)$-equivariant 1-cell. 
Hence there is an $U(1)$-equivariant cofibration sequence 
$$
U(1)_+ \rightarrow S^0 \rightarrow S^c \; \; \eqno{(3.2.6)}
$$ 
where the zero dimensional sphere $S^0$ consists of the north and south pole of the sphere $S^c$. 
Applying the equivariant mapping functor 
$ \Map_{U(1)}^0 ( \underline {\quad } , \, \g {1})$ to the sequence 3.2.6 
gives the fibration sequence 

$$
\Map^0_{U(1)} (U(1)_+ , \, \g {1} ) \leftarrow 
\Map^0_{U(1)} (S^0, \, \g {1} ) 
\leftarrow 
\Map^0_{U(1)} (S^c, \, \g {1} ) \; . \eqno{(3.2.7)} 
$$
 
We have the homeomorphisms 
$$ 
\eqalign 
{\Map^0_{U(1)} (U(1)_+ , \, \g {1} ) &\cong \g {1} =BSU(2) \cr 
\hbox{and} \; \; \; \; \; 
\Map^0_{U(1)} (S^0, \, \g {1} ) &\cong \g {1}^{U(1)} \; . \cr }
$$ 

\proclaim {Lemma 3.2.8} There is a homeomorphism 
$ \g {1}^{U(1)} \cong G_1 (\C ^{\infty}) = BU(1)$. 
\endproclaim

Hence the mapping space 
$\Map^0_{U(1)} (S^c, \, \g {1})$ is the homotopy fiber of the standard map 
$BU(1) \rightarrow BSU(2)$ which is, up to homotopy, the space $S^2$. 
This proves lemma 3.2.5 and finishes the proof of 3.2.3.\enddemo 

\demo{Proof of lemma 3.2.8}  Let $\tau : G_1 (\C^ {\infty})
\hookrightarrow \g {1}$ be the map induced by the natural inclusion $
\C \hookrightarrow \H$.  Recall that we view the space $\H^{\infty}$
as a right quaternionic vector space endowed with the $U(1)$-action
given by left multiplication of weight p. Hence the image $ \tau ( G_1
(\C^ {\infty}))$ lies in the fixed point set $ \g {1}^{U(1)}$ of the
induced circle action on the Grassmannian $\g {1}$ and the map $\tau$
induces a continuous map $\tilde{\tau}: G_1 (\C^{\infty})
\hookrightarrow \g {1}^{U(1)}$. 

Let $w= [x_0 , x_1 , \ldots ]$ be a quaternionic line in $\H^{\infty}$
fixed under the circle action. Without loss of generality we may
assume $x_0 =1$. Given $\lambda \in U(1)$ there is an element $\alpha
(\lambda ) \in Sp(1)$ such that $\lambda^p x_i = x_i \alpha (\lambda
)$ for all $i=0,1, \ldots$. For $i=0$ this gives $\alpha (\lambda )=
\lambda^p $.
Hence $\lambda^p x_i = x_i \lambda^p$ for all $i=0,1, \ldots$ and
hence $x_i$ lies in the centralizer of $U(1)$ in the quaternions $\H$
which is equal to $\C$. 

Hence $x_i \in \C$ for all $i=0,1, \ldots$ and the element $w= [x_0 ,
x_1 , \ldots ]$ lies in the image of the map $\tilde{\tau} : G_1 (\C
^{\infty})\rightarrow \g {1}^{U(1)}$. The map $\tilde{\tau} $ induces
a homeomorphism $G_1 (\C ^{\infty}) \cong
\g {1} ^{U(1)}$. 
This finishes the proof of lemma 3.2.8.\enddemo

\subhead {\rm 3.3} Identification of the component
\endsubhead

We now finish the proof of theorem 3.2.1.  Recall that we denoted by
$\eta = (P \rightarrow S^4)$ a fixed $U(1)$-equivariant principal
$SU(2)$-bundle over the 4-sphere $S^4$ with second Chern class $c_2
(\eta ) =n$. Restricted to the fixed point set $S^2 \subset S^4$ of
the $U(1)$-action on the 4-sphere the bundle $\eta$ is induced from a
$U(1)$-bundle $\nu$ with first Chern class $c_1 (\nu )=k$.  In 3.1.1
we defined the number $k$ to be the charge of the bundle $\eta$.

Let $\Gamma : \Map^0_{U(1)} (S^4 , \, \g {1})
\rightarrow \Omega^2 S^2$ be the map inducing the homotopy equivalence 
in 3.2.3. We want to determine to which component of the double loop space 
$\Omega^2 S^2$ the map $\Gamma$ maps the connected component 
$\Map^0_{U(1)} (S^4 , \, \g {1})^{\eta}$ of the mapping space 
$\Map^0_{U(1)} (S^4 , \, \g {1})$. 
For this it is enough
to compute the degree of the map 
$\Gamma (f): S^2 \rightarrow S^2$ for one map 
$f \in\Map^0_{U(1)} (S^4 , \, \g {1})^{\eta}$. 

We identify the 4-sphere $S^4$ with the smash $S^2 \wedge S^c$ as before. 
Here the circle $U(1)$ acts trivially on the first factor and by 
rotation around the axis through the north-and south pole on the 
second factor. 
The fixed point set of the circle action on $S^4$ is then the subspace 
$\{ [z, S] \in S^2 \wedge S^c \; : z \in S^2 \}$ which is just 
isomorphic to a copy of $S^2$.  
Here $S$ denotes the south  
pole of the second factor $S^c$.  

Denote by $ j: S^2 \rightarrow \Map^0_{U(1)} (S^c , \, \g {1})$ a map
inducing the homotopy equivalence from 3.2.5 and let $\iota^{\ast} :
\Map^0_{U(1)} (S^c , \, \g {1}) \rightarrow
\Map^0_{U(1)} (S^0 , \, \g {1})=BU(1) $ 
the map on the mapping spaces induced by the cofibration 
$ \iota : S^0 \hookrightarrow S^c$ from 3.2.6.   
Let $ \pi = \iota^{\ast} \circ j : S^2 \rightarrow BU(1)$ the composition. 
It is clear that the map $\pi$ induces the identity on second cohomology. 
Hence in order to calculate the degree of the map $\Gamma (f)$ it is enough 
to calculate the map 
$$
(\pi \circ \Gamma (f) )^{\ast} : H^2 (BU(1) , \Z) \rightarrow  
H^2 (S^2 , \Z ) 
\eqno{(3.3.1)}  
$$ 
By construction the map $\pi \circ \Gamma (f)$ is a lift of the 
classifying map $f$ of the bundle 
$\eta$ to the fixed point set $S^2$ 
via the map $B(\iota ) : BU(1) \hookrightarrow BSU(2)$. 
Since the bundle $\eta$ restricted to the fixed point set $S^2 \subset S^4$ 
is induced from a $U(1)$-bundle $\nu$ with first Chern class 
$c_1 ( \nu ) =k $ the map 
$\pi \circ \Gamma (f) : S^2 \rightarrow BU(1)$ is a classifying map for 
the bundle 
$\nu$ and hence the map $(\pi \circ \Gamma (f))^{\ast}$ in 3.3.1 is 
multiplication by $k$. This finally finishes the proof of theorem 3.2.1.  

\Refs 

\item{[A]}Atiyah, M.F.: Magnetic monopoles on hyperbolic spaces. In Proc. 
of Bombay Colloquium 1984 on "Vector Bundles on Algebraic Varieties", 
1987, 1-34, Oxford: Oxford University Press. 

\item{[AB]} Atiyah, M.F., Bott, R.: The Yang-Mills equation over 
Riemann surfaces, Phil. Trans. R. Soc. Lond. A 308 (1982), 523-615. 

\item{[Ba]} Braam, P.J.: Magnetic Monopoles on three-manifolds, 
J. Differential Geometry Vol. 30 (1989), 425-464. 

\item{[DK]}{Donaldson, S.K., Kronheimer, P.B.: ``The Geometry of 
Four-Manifolds'', 
Oxford mathematical monographs, 1990, Oxford: Oxford University Press.} 

\item{[G]} Gritsch, U.: Morse theory for the Yang-Mills functional 
via equivariant homotopy theory, 1997, preprint. 

\item{[Se]}{Segal, G.B.: Equivariant K-theory, Publ. Math. Inst. 
Hautes \' Etudes Sci. 34 
(1968), 129-151.} 

\item{[T]} Taubes, C.H.: 
Monopoles and Maps from $S^2$ to $S^2$, the Topology of 
the Configuration Space, Commun. Math. Phys. 95 (1984), 345-391. 

\item{[Wa]}{Wasserman, A.G.: Equivariant differential topology, Topology 
Vol. 8, 127-150.}  

\endRefs

\enddocument